\documentclass[a4paper]{article}

\usepackage{theorem}
\newtheorem{lem}{Lemma}[section]

\usepackage{amssymb}

\def\GAP{\textsf{GAP}}

\def\PGL{{\rm PGL}}
\def\PSL{{\rm PSL}}

\def\tthdump#1{#1}
\tthdump{}

\def\M{{\cal M}}
\tthdump{\def\MM{{\tilde{\M}}}}

\def\cl{{\rm cl}}
\def\num_known_monster_maxes{$44$}

\begin{document}

\tthdump{\title{{\GAP} Computations Concerning Hamiltonian Cycles in the Generating Graphs of Finite Groups}}

\author{\textsc{Thomas Breuer} \\[0.5cm]
\textit{Lehrstuhl D f\"ur Mathematik} \\
\textit{RWTH Aachen University, 52056 Aachen, Germany}}

\date{April 24th, 2012}

\maketitle

\abstract{
This is a collection of examples showing how
the {\GAP} system~\cite{GAP}
can be used to compute information about the generating graphs of
finite groups.
It includes all examples that were needed for the computational results
in~\cite{GMN}.

The purpose of this writeup is twofold.
On the one hand, the computations are documented this way.
On the other hand, the {\GAP} code shown for the examples can be used as
test input for automatic checking of the data and the functions used.

A first version of this document, which was based on {\GAP}~4.4.12,
had been accessible in the web since November~2009
and is available in the arXiv (no. 0911.5589) since November~2009.
The differences to the current version are as follows.

\begin{itemize}
\item
  The format of the {\GAP} output was adjusted to the changed behaviour
  of {\GAP}~4.5.
\item
  The sporadic simple Monster group has exactly one class of maximal
  subgroups of the type $\PSL(2, 41)$ (see~\cite{NW12}),
  and has no maximal subgroups which have the socle $\PSL(2, 27)$
  (see~\cite{Wil10}).
  As a consequence, the lower bounds computed in Section~\ref{Monster}
  have been improved.
\end{itemize}}

\textwidth16cm
\oddsidemargin0pt

\parskip 1ex plus 0.5ex minus 0.5ex
\parindent0pt

\tableofcontents


\section{Overview}

The purpose of this note is to document the {\GAP} computations
that were carried out in order to obtain the computational results
in~\cite{GMN}.

In order to keep this note self-contained,
we first describe the theory needed, in Section~\ref{background}.
The translation of the relevant formulae into {\GAP} functions
can be found in Section~\ref{functions}.
Then Section~\ref{chartheor} describes the computations that only require
(ordinary) character tables in the
{\GAP} Character Table Library~\cite{CTblLib1.2}.
Computations using also the groups are shown in Section~\ref{grouptheor}.

The examples use the {\GAP} Character Table Library
and the {\GAP} Library of Tables of Marks,
so we first load these packages in the required versions.

Also, we force the assertion level to zero;
this is the default in interactive {\GAP} sessions
but the level is automatically set to $1$
when a file is read with \verb|ReadTest|.

\begin{verbatim}
    gap> if not CompareVersionNumbers( GAPInfo.Version, "4.5" ) then
    >      Error( "need GAP in version at least 4.5" );
    >    fi;
    gap> LoadPackage( "ctbllib", "1.2" );
    true
    gap> LoadPackage( "tomlib", "1.1.1" );
    true
    gap> SetAssertionLevel( 0 );
\end{verbatim}

\section{Theoretical Background}\label{background}

Let $G$ be a finite noncyclic group
and denote by $G^{\times}$ the set of nonidentity elements in $G$.
We define the \emph{generating graph} $\Gamma(G)$ as the undirected graph
on the vertex set $G^{\times}$ by joining two elements $x, y \in G^{\times}$
by an edge if and only if $\langle x, y \rangle = G$ holds.
For $x \in G^{\times}$, the \emph{vertex degree} $d(\Gamma, x)$ is
$|\{ y \in G^{\times}; \langle x, y \rangle = G \}|$.
The \emph{closure} $\cl(\Gamma)$ of the graph $\Gamma$ with $m$ vertices
is defined as the graph with the same vertex set as $\Gamma$,
where the vertices $x, y$ are joined by an edge
if they are joined by an edge in $\Gamma$
or if $d(\Gamma, x) + d(\Gamma, y) \geq m$.
We denote iterated closures by $\cl^{(i)}(\Gamma) = \cl(\cl^{(i-1)}(\Gamma))$,
where $\cl^{(0)}(\Gamma) = \Gamma$.

In the following, we will show that the generating graphs of the following
groups contain a Hamiltonian cycle:
\begin{itemize}
\item
    Nonabelian simple groups of orders at most $10^7$,
\item
    groups $G$ containing a unique minimal normal subgroup $N$
    such that $N$ has order at most $10^6$, $N$ is nonsolvable,
    and $G/N$ is cyclic,
\item
    sporadic simple groups and their automorphism groups.
\end{itemize}

Clearly the condition that $G/N$ is cyclic for all nontrivial normal
subgroups $N$ of $G$ is necessary for $\Gamma(G)$ being connected,
and~\cite[Conjecture~1.6]{GMN} states that this condition is also
sufficient.
By~\cite[Proposition~1.1]{GMN},
this conjecture is true for all solvable groups,
and the second entry in the above list implies that this conjecture holds
for all nonsolvable groups of order up to $10^6$.

The question whether a graph $\Gamma$ contains a Hamiltonian cycle
(i.~e., a closed path in $\Gamma$ that visits each vertex exactly once)
can be answered using the following sufficient criteria (see~\cite{GMN}).
Let $d_1 \leq d_2 \leq \cdots \leq d_m$ be the vertex degrees in $\Gamma$.

\begin{description}
\item[P\'osa's criterion:]
    If $d_k \geq k+1$ holds for $1 \leq k < m/2$
    then $\Gamma$ contains a Hamiltonian cycle.
\item[Chv\'atal's criterion:]
    If $d_k \geq k+1$ or $d_{m-k} \geq m-k$ holds for $1 \leq k < m/2$
    then $\Gamma$ contains a Hamiltonian cycle.
\item[Closure criterion:]
    A graph contains a Hamiltonian cycle if and only if its closure
    contains a Hamiltonian cycle.
\end{description}

%
%

\subsection{Character-Theoretic Lower Bounds for Vertex Degrees}

Using character-theoretic methods similar to those used to obtain the
results in~\cite{BGK}
(the computations for that paper are shown in~\cite{ProbGenArxiv}),
we can compute lower bounds for the vertex degrees in generating graphs,
as follows.

Let $R$ be a set of representatives of conjugacy classes of nonidentity
elements in $G$,
fix $s \in G^{\times}$,
let $\M(G,s)$ denote the set of those maximal subgroups of $G$
that contain $s$,
let $\MM(G,s)$ denote a set of representatives in $\M(G,s)$
w.~r.~t.~conjugacy in $G$.
For a subgroup $M$ of $G$, the \emph{permutation character} $1_M^G$
is defined by
\[
   1_M^G(g):= \frac{|G| \cdot |g^G \cap M|}{|M| \cdot |g^G|},
\]
where $g^G = \{ g^x; x \in G \}$, with $g^x = x^{-1} g x$,
denotes the conjugacy class of $g$ in $G$.
So we have $1_M^G(1) = |G|/|M|$
and thus $|g^G \cap M| = |g^G| \cdot 1_M^G(g) / 1_M^G(1)$.

Doubly counting the set $\{ (s^x, M^y); x, y \in G, s^x \in M^y \}$ yields
$|M^G| \cdot |s^G \cap M| = |s^G| \cdot |\{ M^x; x \in G, s \in M^x \}|$
and thus
$|\{ M^x; x \in G, s \in M^x \}| = |M^G| \cdot 1_M^G(s) / 1_M^G(1) \leq 1_M^G(s)$.
(If $M$ is a \emph{maximal} subgroup of $G$ then either $M$ is normal in $G$
or self-normalizing, and in the latter case the inequality is in fact an
equality.)

Let $\Pi$ denote the multiset of \emph{primitive} permutation characters
of $G$, i.~e.,
of the permutation characters $1_M^G$ where $M$ ranges over representatives
of the conjugacy classes of maximal subgroups of $G$.

Define
\[
   \delta(s, g^G):= |g^G| \cdot \max\left\{ 0, 1 - \sum_{\pi \in \Pi}
                                     \pi(g) \cdot \pi(s) / \pi(1) \right\}
\]
and $d(s, g^G):= |\{ x \in g^G; \langle s, x \rangle = G \}|$,
the contribution of the class $g^G$ to the vertex degree of $s$.
Then we have $d(\Gamma(G), s) = \sum_{x \in R} d(s, x^G)$ and
\begin{eqnarray*}
   d(s, g^G) & = & |g^G| - |\bigcup_{M \in \M(G,s)}
                      \{ x \in g^G; \langle x, s \rangle \subseteq M \}| \\
          & \geq & \max\left\{ 0, |g^G| - \sum_{M \in \M(G,s)}
                                            |g^G \cap M| \right\} \\
             & = & |g^G| \cdot \max\left\{ 0, 1 - \sum_{M \in \M(G,s)}
                                            1_M^G(g) / 1_M^G(1) \right\} \\
          & \geq & |g^G| \cdot \max\left\{ 0, 1 - \sum_{M \in \MM(G,s)}
                             1_M^G(g) \cdot 1_M^G(s) / 1_M^G(1) \right\} \\
             & = & \delta(s, g^G)
\end{eqnarray*}

So $\delta(s):= \sum_{x \in R} \delta(s, x^G)$ is a lower bound
for the vertex degree of $s$; this bound can be computed if $\Pi$ is known.

For computing the vertex degrees of the iterated closures of $\Gamma(G)$,
we define $d^{(0)}(s, g^G):= d(s, g^G)$ and
\[
    d^{(i+1)}(s, g^G):= \left\{ \begin{array}{lcl}
         |g^G|           & ; & d^{(i)}(\Gamma(G), s) + d^{(i)}(\Gamma(G), g)
                               \geq |G|-1 \\
         d^{(i)}(s, g^G) & ; & \mbox{\rm otherwise}
                               \end{array} \right..
\]
Then $d(\cl^{(i)}(\Gamma(G)), s) = \sum_{x \in R} d^{(i)}(s, g^G)$ holds.

Analogously, we set $\delta^{(0)}(s, g^G):= \delta(s, g^G)$,
\[
    \delta^{(i+1)}(s, g^G):= \left\{ \begin{array}{lcl}
         |g^G|          & ; & \delta^{(i)}(s) + \delta^{(i)}(g) \geq |G|-1 \\
         \delta^{(i)}(s, g^G) & ; & \mbox{\rm otherwise}
                             \end{array} \right.
\]
and $\delta^{(i)}(s):= \sum_{x \in R} \delta^{(i)}(s, x^G)$,
a lower bound for $d(\cl^{(i)}(\Gamma(G)), s)$
that can be computed if $\Pi$ is known.

\subsection{Checking the Criteria}\label{critcheck}

Let us assume that we know lower bounds $\beta(s)$ for the vertex degrees
$d(\cl^{(i)}(\Gamma(G)), s)$, for some fixed $i$,
and let us choose representatives $s_1, s_2, \ldots, s_l$ of the nonidentity
conjugacy classes of $G$
such that $\beta(s_1) \leq \beta(s_2) \leq \cdots \leq \beta(s_l)$ holds.
Let $c_k = |s_k^G|$ be the class lengths of these representatives.

Then the first $c_1$ vertex degrees, ordered by increasing size,
are larger than or equal to $\beta(s_1)$,
the next $c_2$ vertex degrees are larger than or equal to $\beta(s_2)$,
and so on.

Then the set of indices in the $k$-th nonidentity class of $G$
for which P\'osa's criterion is not guaranteed by the given bounds is
\[
   \{ x; c_1 + c_2 + \cdots + c_{k-1} < x \leq c_1 + c_2 + \ldots c_k,
         x < (|G| - 1) / 2, \beta(s_k) < x+1 \}.
\]
This is an interval $\{ L_k, L_k + 1, \ldots, U_k \}$ with
\[
   L_k = \max\left\{ 1 + c_1 + c_2 + \cdots + c_{k-1},
                     \beta(s_k)  \right\}
\]
and
\[
   U_k = \min\left\{ c_1 + c_2 + \cdots + c_k,
                     \left\lfloor |G|/2 \right\rfloor - 1 \right\} .
\]
(Note that the generating graph has $m = |G|-1$ vertices,
and that $x < m/2$ is equivalent to
$x \leq \left\lfloor |G|/2 \right\rfloor - 1$.)

The generating graph $\Gamma(G)$ satisfies P\'osa's criterion
if all these intervals are empty,
i.~e., if $L_k > U_k$ holds for $1 \leq k \leq l$.

The set of indices for which Chv\'atal's criterion is not guaranteed
is the intersection of
\[
   \{ m-k; 1 \leq m-k < m/2, d_k < k \}
\]
with the set of indices for which P\'osa's criterion is not guaranteed.

Analogously to the above considerations,
the set of indices $m-x$ in the former set for which Chv\'atal's criterion
is not guaranteed by the given bounds
and such that $x$ is an index in the $k$-th nonidentity class of $G$
is
\[
   \{ m-x; c_1 + c_2 + \cdots + c_{k-1} < x \leq c_1 + c_2 + \ldots c_k,
         1 \leq m-x < (|G| - 1) / 2, \beta(s_k) < x \}.
\]
This is again an interval
$\{ L^{\prime}_k, L^{\prime}_k + 1, \ldots, U^{\prime}_k \}$
with
\[
   L^{\prime}_k = \max\left\{ 1, m - ( c_1 + c_2 + \cdots + c_k ) \right\}
\]
and
\[
   U^{\prime}_k = \min\left\{ m - ( c_1 + c_2 + \cdots + c_{k-1} ) - 1,
                     \left\lfloor |G|/2 \right\rfloor - 1,
                     m-1 - \beta(s_k) \right\} .
\]
The generating graph $\Gamma(G)$ satisfies Chv\'atal's criterion if
the union of the intervals
$\{ L^{\prime}_k, L^{\prime}_k + 1, \ldots, U^{\prime}_k \}$,
for $1 \leq k \leq l$
is disjoint to the union of the intervals $\{ L_k, L_k + 1, \ldots, U_k \}$,
for $1 \leq k \leq l$.

\section{{\GAP} Functions for the Computations}\label{functions}

We describe two approaches to compute, for a given group $G$,
vertex degrees for the generating graph of $G$ or lower bounds for them,
by calculating exact vertex degrees from $G$ itself
(see Section~\ref{groups})
or by deriving lower bounds for the vertex degrees using just
character-theoretic information about $G$ and its subgroups
(see Section~\ref{characters}).
Finally, Section~\ref{clos} deals with deriving lower bounds
of vertex degrees of iterated closures.

\subsection{Computing Vertex Degrees from the Group}\label{groups}

In this section,
the task is to compute the vertex degrees $d(s,g^G)$ using explicit
computations with the group $G$.

The function \verb|IsGeneratorsOfTransPermGroup| checks whether the permutations
in the list \verb|list| generate the permutation group \verb|G|,
\emph{provided that} \verb|G| is transitive on its moved points.
(Note that testing the necessary condition that the elements in \verb|list|
generate a transitive group is usually much faster than testing generation.)
This function has been used already in~\cite{ProbGenArxiv}.

\begin{verbatim}
    gap> BindGlobal( "IsGeneratorsOfTransPermGroup", function( G, list )
    >     local S;
    > 
    >     if not IsTransitive( G ) then
    >       Error( "<G> must be transitive on its moved points" );
    >     fi;
    >     S:= SubgroupNC( G, list );
    > 
    >     return IsTransitive( S, MovedPoints( G ) ) and Size( S ) = Size( G );
    > end );
\end{verbatim}

The function \verb|VertexDegreesGeneratingGraph| takes a
\emph{transitive} permutation group \verb|G|
(in order to be allowed to use \verb|IsGeneratorsOfTransPermGroup|),
the list \verb|classes| of conjugacy classes of \verb|G|
(in order to prescribe an ordering of the classes),
and a list \verb|normalsubgroups| of proper normal subgroups of \verb|G|,
and returns the matrix $[ d(s, g^{\verb|G|}) ]_{s, g}$ of vertex degrees,
with rows and columns indexed by nonidentity class representatives
ordered as in the list \verb|classes|.
(The class containing the identity element may be contained in \verb|classes|.)

The following criteria are used in this function.
\begin{itemize}
\item
    The function tests the (non)generation only for representatives of
    $C_G(g)$-$C_G(s)$-double cosets,
    where $C_G(g):= \{ x \in G; g x = x g \}$
    denotes the centralizer of $g$ in $G$.
    Note that for $c_1 \in C_G(g)$, $c_2 \in C_G(s)$,
    and a representative $r \in G$, we have
    $\langle s, g^{c_1 r c_2} \rangle = \langle s, g^r \rangle^{c_2}$.
    If $\langle s, g^r \rangle = G$ then
    the double coset $D = C_G(g) r C_G(s)$ contributes $|D|/|C_G(g)|$
    to the vertex degree $d(s, g^G)$,
    otherwise the contribution is zero.
\item
    We have $d(s, g^G) \cdot |C_G(g)| = d(g, s^G) \cdot |C_G(s)|$.
    (To see this,
    either establish a bijection of the above double cosets,
    or doubly count the edges between elements of the conjugacy classes
    of $s$ and $g$.)
\item
    If $\langle s_1 \rangle = \langle s_2 \rangle$ and
    $\langle g_1 \rangle = \langle g_2 \rangle$ hold then we have
    $d(s_1, g_1^G) = d(s_2, g_1^G) = d(s_1, g_2^G) = d(s_2, g_2^G)$,
    so only one of these values must be computed.
\item
    If both $s$ and $g$ are contained in one of the normal subgroups given
    then $d(s, g^G)$ is zero.
\item
    If $G$ is not a dihedral group and both $s$ and $g$ are involutions
    then $d(s, g^G)$ is zero.
\end{itemize}

\begin{verbatim}
    gap> BindGlobal( "VertexDegreesGeneratingGraph",
    >     function( G, classes, normalsubgroups )
    >     local nccl, matrix, cents, powers, normalsubgroupspos, i, j, g_i, nsg,
    >           g_j, gen, pair, d, pow;
    > 
    >     if not IsTransitive( G ) then
    >       Error( "<G> must be transitive on its moved points" );
    >     fi;
    > 
    >     classes:= Filtered( classes, C -> Order( Representative( C ) ) <> 1 );
    >     nccl:= Length( classes );
    >     matrix:= [];
    >     cents:= [];
    >     powers:= [];
    >     normalsubgroupspos:= [];
    >     for i in [ 1 .. nccl ] do
    >       matrix[i]:= [];
    >       if IsBound( powers[i] ) then
    >         # The i-th row equals the earlier row `powers[i]'.
    >         for j in [ 1 .. i ] do
    >           matrix[i][j]:= matrix[ powers[i] ][j];
    >           matrix[j][i]:= matrix[j][ powers[i] ];
    >         od;
    >       else
    >         # We have to compute the values.
    >         g_i:= Representative( classes[i] );
    >         nsg:= Filtered( [ 1 .. Length( normalsubgroups ) ],
    >                         i -> g_i in normalsubgroups[i] );
    >         normalsubgroupspos[i]:= nsg;
    >         cents[i]:= Centralizer( G, g_i );
    >         for j in [ 1 .. i ] do
    >           g_j:= Representative( classes[j] );
    >           if IsBound( powers[j] ) then
    >             matrix[i][j]:= matrix[i][ powers[j] ];
    >             matrix[j][i]:= matrix[ powers[j] ][i];
    >           elif not IsEmpty( Intersection( nsg, normalsubgroupspos[j] ) )
    >                or ( Order( g_i ) = 2 and Order( g_j ) = 2
    >                     and not IsDihedralGroup( G ) ) then
    >             matrix[i][j]:= 0;
    >             matrix[j][i]:= 0;
    >           else
    >             # Compute $d(g_i, g_j^G)$.
    >             gen:= 0;
    >             for pair in DoubleCosetRepsAndSizes( G, cents[j], cents[i] ) do
    >               if IsGeneratorsOfTransPermGroup( G, [ g_i, g_j^pair[1] ] ) then
    >                 gen:= gen + pair[2];
    >               fi;
    >             od;
    >             matrix[i][j]:= gen / Size( cents[j] );
    >             if i <> j then
    >               matrix[j][i]:= gen / Size( cents[i] );
    >             fi;
    >           fi;
    >         od;
    > 
    >         # For later, provide information about algebraic conjugacy.
    >         for d in Difference( PrimeResidues( Order( g_i ) ), [ 1 ] ) do
    >           pow:= g_i^d;
    >           for j in [ i+1 .. nccl ] do
    >             if not IsBound( powers[j] ) and pow in classes[j] then
    >               powers[j]:= i;
    >               break;
    >             fi;
    >           od;
    >         od;
    >       fi;
    >     od;
    > 
    >     return matrix;
    > end );
\end{verbatim}

\subsection{Computing Lower Bounds for Vertex Degrees}\label{characters}

In this section,
the task is to compute the lower bounds $\delta(s, g^G)$ for the
vertex degrees $d(s, g^G)$ using character-theoretic methods.

We provide {\GAP} functions for computing
the multiset $\Pi$ of the primitive permutation characters
of a given group $G$ and
for computing the lower bounds $\delta(s, g^G)$ from $\Pi$.

For many almost simple groups, the {\GAP} libraries of character tables
and of tables of marks contain information for quickly computing
the primitive permutation characters of the group in question.
Therefore, the function \verb|PrimitivePermutationCharacters| takes as its
argument not the group $G$ but its character table $T$, say.
(This function is shown already in~\cite{ProbGenArxiv}.)

If $T$ is contained in the {\GAP} Character Table Library
(see~\cite{CTblLib1.2})
then the complete set of primitive permutation characters
can be easily computed if the character tables of all maximal
subgroups and their class fusions into $T$ are known
(in this case, we check whether the attribute \verb|Maxes| of $T$ is bound)
or if the table of marks of $G$ and the class fusion from $T$ into this
table of marks are known
(in this case, we check whether the attribute \verb|FusionToTom| of $T$ is bound).
If the attribute \verb|UnderlyingGroup| of $T$ is bound then the group
stored as the value of this attribute
can be used to compute the primitive permutation characters.
The latter happens if $T$ was computed from the group $G$;
for tables in the {\GAP} Character Table Library,
this is not the case by default.

The {\GAP} function \verb|PrimitivePermutationCharacters| tries to compute
the primitive permutation characters of a group using this information;
it returns the required list of characters if this can be computed this way,
otherwise \verb|fail| is returned.
(For convenience, we use the {\GAP} mechanism of \emph{attributes}
in order to store the permutation characters in the character table object
once they have been computed.)

\begin{verbatim}
    gap> DeclareAttribute( "PrimitivePermutationCharacters", IsCharacterTable );
    gap> InstallMethod( PrimitivePermutationCharacters,
    >     [ IsCharacterTable ],
    >     function( tbl )
    >     local maxes, i, fus, poss, tom, G;
    > 
    >     if HasMaxes( tbl ) then
    >       maxes:= List( Maxes( tbl ), CharacterTable );
    >       for i in [ 1 .. Length( maxes ) ] do
    >         fus:= GetFusionMap( maxes[i], tbl );
    >         if fus = fail then
    >           fus:= PossibleClassFusions( maxes[i], tbl );
    >           poss:= Set( List( fus, map -> InducedClassFunctionsByFusionMap(
    >              maxes[i], tbl, [ TrivialCharacter( maxes[i] ) ], map )[1] ) );
    >           if Length( poss ) = 1 then
    >             maxes[i]:= poss[1];
    >           else
    >             return fail;
    >           fi;
    >         else
    >           maxes[i]:= TrivialCharacter( maxes[i] )^tbl;
    >         fi;
    >       od;
    >       return maxes;
    >     elif HasFusionToTom( tbl ) then
    >       tom:= TableOfMarks( tbl );
    >       maxes:= MaximalSubgroupsTom( tom );
    >       return PermCharsTom( tbl, tom ){ maxes[1] };
    >     elif HasUnderlyingGroup( tbl ) then
    >       G:= UnderlyingGroup( tbl );
    >       return List( MaximalSubgroupClassReps( G ),
    >                    M -> TrivialCharacter( M )^tbl );
    >     fi;
    > 
    >     return fail;
    > end );
\end{verbatim}


The next function computes the lower bounds $\delta(s, g^G)$ from
the two lists \verb|classlengths| of conjugacy class lengths of the group $G$
and \verb|prim| of all primitive permutation characters of $G$.
(The first entry in \verb|classlengths| is assumed to represent the class
containing the identity element of $G$.)
The return value is the matrix that contains in row $i$ and column $j$
the value $\delta(s, g^G)$, where $s$ and $g$ are in the conjugacy classes
represented by the $(i+1)$-st and $(j+1)$-st column of \verb|tbl|, respectively.
So the row sums of this matrix are the values $\delta(s)$.

\begin{verbatim}
    gap> LowerBoundsVertexDegrees:= function( classlengths, prim )
    >     local sizes, nccl;
    > 
    >     nccl:= Length( classlengths );
    >     return List( [ 2 .. nccl ],
    >              i -> List( [ 2 .. nccl ],
    >                     j -> Maximum( 0, classlengths[j] - Sum( prim,
    >                     pi -> classlengths[j] * pi[j] * pi[i] / pi[1] ) ) ) );
    > end;;
\end{verbatim}

\subsection{Evaluating the (Lower Bounds for the) Vertex Degrees}\label{clos}

In this section,
the task is to compute (lower bounds for) the vertex degrees of iterated
closures of a generating graph from (lower bounds for) the vertex degrees
of the graph itself,
and then to check the criteria of P\'osa and Chv\'atal.

The arguments of all functions defined in this section are
the list \verb|classlengths| of conjugacy class lengths for the group $G$
(including the class of the identity element, in the first position)
and a matrix \verb|bounds| of the values $d^{(i)}(s, g^G)$
or $\delta^{(i)}(s, g^G)$,
with rows and columns indexed by nonidentity class representatives
$s$ and $g$, respectively.
Such a matrix is returned by the functions \verb|VertexDegreesGeneratingGraph| or
\verb|LowerBoundsVertexDegrees|, respectively.

The function \verb|LowerBoundsVertexDegreesOfClosure| returns the corresponding
matrix of the values $d^{(i+1)}(s, g^G)$ or $\delta^{(i+1)}(s, g^G)$,
respectively.

\begin{verbatim}
    gap> LowerBoundsVertexDegreesOfClosure:= function( classlengths, bounds )
    >     local delta, newbounds, size, i, j;
    > 
    >     delta:= List( bounds, Sum );
    >     newbounds:= List( bounds, ShallowCopy );
    >     size:= Sum( classlengths );
    >     for i in [ 1 .. Length( bounds ) ] do
    >       for j in [ 1 .. Length( bounds ) ] do
    >         if delta[i] + delta[j] >= size - 1 then
    >           newbounds[i][j]:= classlengths[ j+1 ];
    >         fi;
    >       od;
    >     od;
    > 
    >     return newbounds;
    > end;;
\end{verbatim}

Once the values $d^{(i)}(s, g^G)$ or $\delta^{(i)}(s, g^G)$ are known,
we can check whether P\'osa's or Chv\'atal's criterion is satisfied
for the graph $\cl^{(i)}(\Gamma(G))$,
using the function \verb|CheckCriteriaOfPosaAndChvatal| shown below.
(Of course a \emph{negative} result is meaningless in the case that only
lower bounds for the vertex degrees are used.)

The idea is to compute the row sums of the given matrix,
and to compute the intervals $\{ L_k, L_k + 1, \ldots, U_k \}$
and $\{ L^{\prime}_k, L^{\prime}_k + 1, \ldots, U^{\prime}_k \}$
that were introduced in Section~\ref{critcheck}.

The function \verb|CheckCriteriaOfPosaAndChvatal| returns,
given the list of class lengths of $G$ and the matrix of (bounds for the)
vertex degrees, a record with the components
\verb|badForPosa| (the list of those pairs $[ L_k, U_k ]$ with the property
$L_k \leq U_k$),
\verb|badForChvatal| (the list of pairs of lower and upper bounds of
nonempty intervals where Chv\'atal's criterion may be violated),
and \verb|data| (the sorted list of triples $[ \delta(g_k), |g_k^G|, \iota(k) ]$,
where $\iota(k)$ is the row and column position of $g_k$ in the matrix
\verb|bounds|).
The ordering of class lengths must of course be compatible with the
ordering of rows and columns of the matrix,
and the identity element of $G$ must belong to the first entry in the
list of class lengths.

\begin{verbatim}
    gap> CheckCriteriaOfPosaAndChvatal:= function( classlengths, bounds )
    >     local size, degs, addinterval, badForPosa, badForChvatal1, pos, half,
    >           i, low1, upp2, upp1, low2, badForChvatal, interval1, interval2;
    > 
    >     size:= Sum( classlengths );
    >     degs:= List( [ 2 .. Length( classlengths ) ],
    >                  i -> [ Sum( bounds[ i-1 ] ), classlengths[i], i ] );
    >     Sort( degs );
    > 
    >     addinterval:= function( intervals, low, upp )
    >       if low <= upp then
    >         Add( intervals, [ low, upp ] );
    >       fi;
    >     end;
    > 
    >     badForPosa:= [];
    >     badForChvatal1:= [];
    >     pos:= 1;
    >     half:= Int( size / 2 ) - 1;
    >     for i in [ 1 .. Length( degs ) ] do
    >       # We have pos = c_1 + c_2 + \cdots + c_{i-1} + 1
    >       low1:= Maximum( pos, degs[i][1] );  # L_i
    >       upp2:= Minimum( half, size-1-pos, size-1-degs[i][1] ); # U'_i
    >       pos:= pos + degs[i][2];
    >       upp1:= Minimum( half, pos-1 ); # U_i
    >       low2:= Maximum( 1, size-pos ); # L'_i
    >       addinterval( badForPosa, low1, upp1 );
    >       addinterval( badForChvatal1, low2, upp2 );
    >     od;
    > 
    >     # Intersect intervals.
    >     badForChvatal:= [];
    >     for interval1 in badForPosa do
    >       for interval2 in badForChvatal1 do
    >         addinterval( badForChvatal, Maximum( interval1[1], interval2[1] ),
    >                                     Minimum( interval1[2], interval2[2] ) );
    >       od;
    >     od;
    > 
    >     return rec( badForPosa:= badForPosa,
    >                 badForChvatal:= Set( badForChvatal ),
    >                 data:= degs );
    > end;;
\end{verbatim}

Finally, the function \verb|HamiltonianCycleInfo| assumes that the matrix \verb|bounds|
contains lower bounds for the vertex degrees in the generating graph
$\Gamma$, and returns a string that describes the minimal $i$ with the
property that the given bounds suffice to show that $cl^{(i)}(\Gamma)$
satisfies P\'osa's or Chv\'atal's criterion,
if such a closure exists.
If no closure has this property, the string \verb|"no decision"| is returned.

\begin{verbatim}
    gap> HamiltonianCycleInfo:= function( classlengths, bounds )
    >     local i, result, res, oldbounds;
    > 
    >     i:= 0;
    >     result:= rec( Posa:= fail, Chvatal:= fail );
    >     repeat
    >       res:= CheckCriteriaOfPosaAndChvatal( classlengths, bounds );
    >       if result.Posa = fail and IsEmpty( res.badForPosa ) then
    >         result.Posa:= i;
    >       fi;
    >       if result.Chvatal = fail and IsEmpty( res.badForChvatal ) then
    >         result.Chvatal:= i;
    >       fi;
    >       i:= i+1;
    >       oldbounds:= bounds;
    >       bounds:= LowerBoundsVertexDegreesOfClosure( classlengths, bounds );
    >     until oldbounds = bounds;
    > 
    >     if result.Posa <> fail then
    >       if result.Posa <> result.Chvatal then
    >         return Concatenation( "Chvatal for ", Ordinal( result.Chvatal ),
    >             " closure, Posa for ", Ordinal( result.Posa ), " closure" );
    >       else
    >         return Concatenation( "Posa for ", Ordinal( result.Posa ),
    >             " closure" );
    >       fi;
    >     elif result.Chvatal <> fail then
    >       return Concatenation( "Chvatal for ", Ordinal( result.Chvatal ),
    >                             " closure" );
    >     else
    >       return "no decision";
    >     fi;
    > end;;
\end{verbatim}

\section{Character-Theoretic Computations}\label{chartheor}

In this section, we apply the functions introduced in Section~\ref{functions}
to character tables of almost simple groups
that are available in the {\GAP} Character Table Library.

Our first examples are the sporadic simple groups, in Section~\ref{spor},
then their automorphism groups are considered in Section~\ref{sporaut}.
Small alternating and symmetric groups are treated in
Section~\ref{symmalt}.

For our convenience, we provide a small function that takes as its
argument only the character table in question,
and returns a string, either \verb|"no prim. perm. characters"| or the
return value of \verb|HamiltonianCycleInfo| for the bounds computed from
the primitive permutation characters.

\begin{verbatim}
    gap> HamiltonianCycleInfoFromCharacterTable:= function( tbl )
    >     local prim, classlengths, bounds;
    > 
    >     prim:= PrimitivePermutationCharacters( tbl );
    >     if prim = fail then
    >       return "no prim. perm. characters";
    >     fi;
    >     classlengths:= SizesConjugacyClasses( tbl );
    >     bounds:= LowerBoundsVertexDegrees( classlengths, prim );
    >     return HamiltonianCycleInfo( classlengths, bounds );
    > end;;
\end{verbatim}

\subsection{Sporadic Simple Groups, except the Monster}\label{spor}

The {\GAP} Character Table Library contains the tables of maximal subgroups
of all sporadic simple groups except $M$.


So the function \verb|PrimitivePermutationCharacters|
can be used to compute all primitive permutation characters
for $25$ of the $26$ sporadic simple groups.

\begin{verbatim}
    gap> spornames:= AllCharacterTableNames( IsSporadicSimple, true,
    >                    IsDuplicateTable, false );
    [ "B", "Co1", "Co2", "Co3", "F3+", "Fi22", "Fi23", "HN", "HS", "He", "J1", 
      "J2", "J3", "J4", "Ly", "M", "M11", "M12", "M22", "M23", "M24", "McL", 
      "ON", "Ru", "Suz", "Th" ]
    gap> for tbl in List( spornames, CharacterTable ) do
    >      info:= HamiltonianCycleInfoFromCharacterTable( tbl );
    >      if info <> "Posa for 0th closure" then
    >        Print( Identifier( tbl ), ": ", info, "\n" );
    >      fi;
    >    od;
    M: no prim. perm. characters
\end{verbatim}

It turns out that only for the Monster group,
the information available in the {\GAP} Character Table Library
is not sufficient to prove that the generating graph
contains a Hamiltonian cycle.

\subsection{The Monster}\label{Monster}

Currently {\num_known_monster_maxes} classes of maximal subgroups of the
Monster group $M$ are
known, but there may be more, see~\cite{NW12}.
For some of the known ones, the character table is not known,
and for some of those with known character table,
the permutation character is not uniquely determined by the character tables
involved.

Nevertheless, we will show that the generating graph of $M$ satisfies
P\'osa's criterion.
For that, we use the information that is available.

For some of the known maximal subgroups $S$,
the character tables are available in the {\GAP} Character Table Library,
and we derive upper bounds for the values of the primitive
permutation characters $1_S^M$ from the possible class fusions
from $S$ into $M$.
For the other subgroups $S$, the permutation characters $1_S^M$ have been
computed with other methods.

The list \verb|prim| defined below has length {\num_known_monster_maxes}.
The entry at position $i$ is a list of length one or two.
If \verb|prim[|$i$\verb|]| has length one then its unique entry is the identifier
of the library character table of the $i$-th maximal subgroup of $M$.
If \verb|prim[|$i$\verb|]| has length two then its entries are a string describing
the structure of the $i$-th maximal subgroup $S$ of $M$ and the
permutation character $1_S^M$.

(The construction of the explicitly given characters in this list
will be documented elsewhere.
Some of the constructions can be found in~\cite{ctblpope}.)

\begin{verbatim}
    gap> m:= CharacterTable( "M" );;
    gap> primdata:= [
    > [ "2.B" ],
    > [ "2^1+24.Co1" ],
    > [ "3.F3+.2" ],
    > [ "2^2.2E6(2).3.2" ],
    > [ "2^(10+16).O10+(2)",
    >   Character( m, [ 512372707698741056749515292734375,
    >     405589064025344574375, 29628786742129575, 658201521662685,
    >     215448838605, 0, 121971774375, 28098354375, 335229607, 108472455,
    >     164587500, 4921875, 2487507165, 2567565, 26157789, 6593805, 398925, 0,
    >     437325, 0, 44983, 234399, 90675, 21391, 41111, 12915, 6561, 6561,
    >     177100, 7660, 6875, 315, 275, 0, 113373, 17901, 57213, 0, 4957, 1197,
    >     909, 301, 397, 0, 0, 0, 3885, 525, 0, 2835, 90, 45, 0, 103, 67, 43, 28,
    >     81, 189, 9, 9, 9, 0, 540, 300, 175, 20, 15, 7, 420, 0, 0, 0, 0, 0, 0,
    >     0, 165, 61, 37, 37, 0, 9, 9, 13, 5, 0, 0, 0, 0, 0, 0, 77, 45, 13, 0, 0,
    >     45, 115, 19, 10, 0, 5, 5, 9, 9, 1, 1, 0, 0, 4, 0, 0, 9, 9, 3, 1, 0, 0,
    >     0, 0, 0, 0, 4, 1, 1, 0, 24, 0, 0, 0, 0, 0, 6, 0, 0, 0, 0, 0, 0, 1, 0,
    >     4, 0, 0, 0, 0, 1, 0, 0, 0, 0, 0, 3, 3, 1, 1, 2, 0, 3, 3, 0, 0, 0, 0, 0,
    >     0, 0, 0, 0, 0, 0, 0, 2, 0, 0, 0, 0, 0, 0, 0, 0, 0, 0, 0, 0, 0, 0, 0, 0,
    >     0, 0, 0, 0 ] ) ],
    > [ "2^(2+11+22).(M24xS3)",
    >   Character( m, [ 16009115629875684006343550944921875,
    >     7774182899642733721875, 120168544413337875, 4436049512692980,
    >     215448838605, 131873639625, 760550656275, 110042727795, 943894035,
    >     568854195, 1851609375, 0, 4680311220, 405405, 78624756, 14467005,
    >     178605, 248265, 874650, 0, 76995, 591163, 224055, 34955, 29539, 20727,
    >     0, 0, 375375, 15775, 0, 0, 0, 495, 116532, 3645, 62316, 1017, 11268,
    >     357, 1701, 45, 117, 705, 0, 0, 4410, 1498, 0, 3780, 810, 0, 0, 83, 135,
    >     31, 0, 0, 0, 0, 0, 0, 0, 255, 195, 0, 215, 0, 0, 210, 0, 42, 0, 35, 15,
    >     1, 1, 160, 48, 9, 92, 25, 9, 9, 5, 1, 21, 0, 0, 0, 0, 0, 98, 74, 42, 0,
    >     0, 0, 120, 76, 10, 0, 0, 0, 0, 0, 1, 1, 0, 6, 0, 0, 0, 0, 0, 0, 0, 0,
    >     0, 0, 0, 0, 5, 3, 0, 0, 0, 18, 0, 10, 0, 3, 3, 0, 1, 1, 1, 1, 0, 0, 2,
    >     0, 0, 0, 0, 0, 0, 2, 0, 0, 0, 0, 0, 6, 12, 0, 0, 2, 0, 0, 0, 2, 0, 0,
    >     1, 1, 0, 0, 0, 0, 0, 0, 0, 2, 0, 2, 0, 0, 1, 1, 1, 1, 0, 0, 0, 0, 0, 0,
    >     0, 0, 0, 0, 0, 0 ] ) ],
    > [ "3^(1+12).2.Suz.2" ],
    > [ "2^(5+10+20).(S3xL5(2))",
    >   Character( m, [ 391965121389536908413379198941796875,
    >     23914487292951376996875, 474163138042468875, 9500455925885925,
    >     646346515815, 334363486275, 954161764875, 147339103275, 1481392395,
    >     1313281515, 0, 8203125, 9827885925, 1216215, 91556325, 9388791, 115911,
    >     587331, 874650, 0, 79515, 581955, 336375, 104371, 62331, 36855, 0, 0,
    >     0, 0, 28125, 525, 1125, 0, 188325, 16767, 88965, 2403, 9477, 1155, 891,
    >     207, 351, 627, 0, 0, 4410, 1498, 0, 0, 0, 30, 150, 91, 151, 31, 0, 0,
    >     0, 0, 0, 0, 0, 0, 0, 125, 0, 5, 5, 210, 0, 42, 0, 0, 0, 0, 0, 141, 45,
    >     27, 61, 27, 9, 9, 7, 3, 15, 0, 0, 0, 0, 0, 98, 74, 42, 0, 0, 30, 0, 0,
    >     0, 6, 6, 6, 3, 3, 1, 1, 0, 0, 0, 0, 0, 0, 0, 0, 0, 0, 0, 0, 0, 0, 0, 0,
    >     1, 1, 0, 18, 0, 10, 0, 0, 0, 0, 0, 0, 0, 0, 0, 0, 1, 0, 0, 0, 0, 0, 0,
    >     2, 0, 0, 0, 0, 0, 0, 0, 2, 2, 0, 2, 3, 3, 0, 0, 0, 0, 0, 0, 0, 0, 0, 0,
    >     0, 0, 2, 0, 2, 0, 0, 0, 0, 0, 0, 3, 3, 0, 0, 0, 0, 0, 0, 0, 0, 0, 0 ] ) ],
    > [ "S3xTh" ],
    > [ "[2^39].(L3(2)x3S6)",
    >   Character( m, [ 4050306254358548053604918389065234375,
    >     148844831270071996434375, 2815847622206994375, 14567365753025085,
    >     3447181417680, 659368198125, 3520153823175, 548464353255, 5706077895,
    >     3056566695, 264515625, 0, 19572895485, 6486480, 186109245, 61410960,
    >     758160, 688365, 58310, 0, 172503, 1264351, 376155, 137935, 99127,
    >     52731, 0, 0, 119625, 3625, 0, 0, 0, 0, 402813, 29160, 185301, 2781,
    >     21069, 1932, 4212, 360, 576, 1125, 0, 0, 1302, 294, 0, 2160, 810, 0, 0,
    >     111, 179, 43, 0, 0, 0, 0, 0, 0, 0, 185, 105, 0, 65, 0, 0, 224, 0, 14,
    >     0, 0, 0, 0, 0, 337, 105, 36, 157, 37, 18, 18, 16, 4, 21, 0, 0, 0, 0, 0,
    >     70, 38, 14, 0, 0, 0, 60, 40, 10, 0, 0, 0, 0, 0, 1, 1, 0, 0, 0, 10, 0,
    >     0, 0, 0, 0, 0, 0, 0, 0, 0, 5, 1, 0, 0, 0, 24, 0, 6, 0, 0, 0, 0, 0, 0,
    >     0, 0, 0, 0, 3, 0, 0, 0, 0, 0, 0, 2, 0, 0, 0, 0, 0, 6, 8, 0, 0, 2, 0, 0,
    >     0, 0, 0, 0, 0, 0, 2, 0, 0, 0, 0, 0, 0, 4, 0, 2, 0, 0, 0, 0, 0, 0, 0, 0,
    >     0, 0, 0, 0, 0, 0, 4, 0, 0, 0 ] ) ],
    > [ "3^8.O8-(3).2_3",
    >   Character( m, [ 6065553341050124859256025907200000000,
    >     117457246944126566400000, 2373192769339392000, 4237313863946240,
    >     1528370386400, 480247040000, 1485622476800, 207447654400, 3523215360,
    >     1124597760, 2926000000, 0, 4720235520, 18601440, 49864704, 14602080,
    >     1914720, 645120, 0, 168070, 0, 811008, 133120, 204800, 0, 8192, 3521,
    >     4250, 308000, 28800, 0, 0, 0, 0, 53504, 1520, 68992, 3584, 2304, 672,
    >     7216, 240, 192, 960, 156, 0, 0, 0, 70, 7840, 550, 0, 0, 0, 0, 0, 0, 48,
    >     93, 57, 18, 24, 0, 160, 200, 0, 320, 0, 0, 0, 49, 0, 4, 0, 0, 0, 0,
    >     144, 0, 0, 80, 0, 20, 64, 20, 0, 8, 0, 12, 0, 2, 2, 0, 0, 0, 6, 0, 0,
    >     20, 24, 30, 0, 0, 0, 0, 0, 0, 0, 0, 0, 0, 0, 0, 5, 2, 7, 0, 0, 12, 0,
    >     0, 0, 20, 8, 0, 0, 5, 0, 4, 0, 1, 0, 0, 1, 0, 0, 0, 0, 0, 0, 0, 0, 0,
    >     12, 0, 0, 0, 0, 2, 2, 0, 0, 0, 2, 4, 0, 0, 2, 0, 0, 0, 0, 0, 0, 0, 0,
    >     0, 0, 0, 0, 0, 0, 0, 0, 0, 0, 0, 0, 0, 0, 0, 0, 0, 0, 0, 0, 0, 0, 0, 0,
    >     0, 0, 0, 0 ] ) ],
    > [ "(D10xHN).2" ],
    > [ "(3^2:2xO8+(3)).S4",
    >   Character( m,
    >   [ 377694424605514962329798663208960000000, 2359667150587942666240000,
    >     28756421759729664000, 377826645416419584, 16593035298840,
    >     5193792576000, 14007297638400, 1715997638656, 5830082560, 2683699200,
    >     5266800000, 0, 47782831360, 241801560, 626008320, 48633880, 9541080,
    >     483840, 2332400, 0, 16384, 3964928, 926720, 102400, 16384, 32256,
    >     51030, 7371, 800800, 41600, 0, 0, 0, 0, 248640, 120480, 200656, 13440,
    >     13696, 1260, 4708, 1120, 1864, 0, 26, 0, 7840, 336, 0, 4284, 180, 0, 0,
    >     0, 0, 0, 0, 265, 418, 270, 99, 81, 0, 480, 456, 0, 0, 0, 0, 455, 0, 56,
    >     0, 0, 0, 0, 0, 680, 64, 4, 16, 32, 18, 26, 16, 4, 0, 0, 10, 0, 9, 0,
    >     28, 24, 8, 0, 0, 0, 160, 20, 20, 0, 0, 0, 0, 0, 0, 0, 0, 0, 0, 0, 0,
    >     18, 3, 4, 1, 0, 8, 2, 0, 0, 20, 8, 0, 0, 0, 19, 0, 0, 0, 0, 0, 0, 0, 0,
    >     0, 0, 0, 0, 0, 0, 0, 2, 0, 1, 0, 4, 0, 0, 0, 0, 0, 6, 0, 0, 0, 0, 0, 0,
    >     0, 0, 0, 0, 0, 0, 0, 0, 0, 0, 4, 0, 0, 1, 0, 0, 0, 0, 0, 0, 0, 0, 0, 0,
    >     0, 0, 0, 0, 2, 2, 0, 0, 0, 0 ] ) ],
    > [ "3^(2+5+10).(M11x2S4)",
    >   Character( m,
    >   [ 16458603283969466072643078298009600000000, 20359256136981938176000000,
    >     145987312780574720000, 724314893457326080, 21414300718720,
    >     18249387520000, 540226355200, 1703254425600, 4697620480, 4771020800,
    >     23408000000, 0, 43256012800, 98483840, 909246464, 213623680, 8362880,
    >     4444160, 0, 0, 0, 475136, 998400, 81920, 0, 35840, 25312, 10597, 0,
    >     128000, 0, 0, 0, 440, 93184, 160, 134400, 1792, 7168, 560, 15888, 160,
    >     64, 320, 0, 0, 0, 0, 0, 19880, 2240, 0, 0, 0, 0, 0, 0, 301, 148, 200,
    >     221, 53, 0, 640, 0, 0, 0, 0, 0, 0, 0, 0, 0, 0, 40, 0, 0, 224, 0, 0, 32,
    >     32, 20, 156, 8, 0, 8, 0, 0, 0, 1, 1, 0, 0, 0, 0, 0, 0, 0, 104, 80, 0,
    >     0, 0, 0, 0, 0, 0, 8, 8, 0, 0, 0, 4, 1, 6, 1, 0, 0, 0, 0, 0, 0, 16, 0,
    >     0, 0, 0, 0, 0, 0, 4, 4, 2, 0, 0, 0, 0, 0, 0, 0, 0, 0, 0, 0, 1, 0, 0, 0,
    >     0, 0, 0, 0, 0, 4, 0, 0, 0, 0, 0, 0, 0, 4, 0, 0, 0, 0, 0, 0, 0, 0, 0, 0,
    >     0, 0, 0, 0, 0, 2, 2, 0, 0, 0, 0, 0, 0, 0, 0, 0, 0, 0, 0, 0, 0 ] ) ],
    > [ "3^(3+2+6+6):(L3(3)xSD16)",
    >   Character( m,
    >   [ 69632552355255433384259177414656000000000, 10962676381451812864000000,
    >     276872489756262400000, 816070626832384000, 52168426710400,
    >     4994569216000, 29712449536000, 917136998400, 32883343360, 14313062400,
    >     0, 0, 53947801600, 445244800, 995491840, 268777600, 8579200, 2007040,
    >     0, 0, 0, 4505600, 588800, 245760, 0, 35840, 24760, 4105, 0, 0, 0, 0, 0,
    >     0, 148480, 8800, 134400, 1792, 13312, 1680, 20784, 1120, 448, 960, 156,
    >     0, 0, 0, 0, 0, 0, 0, 0, 0, 0, 0, 0, 241, 304, 184, 121, 49, 0, 0, 0, 0,
    >     0, 0, 0, 0, 0, 0, 0, 0, 0, 0, 0, 416, 0, 0, 96, 32, 20, 92, 24, 0, 8,
    >     0, 28, 0, 1, 1, 0, 0, 0, 0, 0, 0, 0, 0, 0, 0, 0, 0, 0, 0, 0, 0, 0, 0,
    >     0, 0, 0, 4, 1, 12, 1, 0, 12, 0, 0, 0, 0, 0, 0, 0, 0, 0, 0, 0, 0, 0, 0,
    >     0, 0, 0, 0, 0, 0, 0, 0, 0, 0, 12, 0, 1, 0, 0, 0, 0, 0, 0, 0, 0, 0, 0,
    >     0, 0, 0, 0, 0, 0, 0, 0, 0, 0, 0, 0, 0, 0, 4, 0, 0, 0, 0, 0, 0, 0, 0, 0,
    >     0, 0, 0, 0, 0, 0, 0, 0, 4, 4, 0, 0, 0, 0 ] ) ],
    > [ "5^(1+6):2.J2.4" ],
    > [ "(7:3xHe):2" ],
    > [ "(A5xA12):2" ],
    > [ "5^(3+3).(2xL3(5))" ],
    > [ "(A6xA6xA6).(2xS4)" ],
    > [ "(A5xU3(8):3):2" ],
    > [ "5^(2+2+4):(S3xGL2(5))" ],
    > [ "(L3(2)xS4(4):2).2" ],
    > [ "7^(1+4):(3x2.S7)" ],
    > [ "(5^2:[2^4]xU3(5)).S3" ],
    > [ "(L2(11)xM12):2" ],
    > [ "(A7x(A5xA5):2^2):2" ],
    > [ "5^4:(3x2.L2(25)).2" ],
    > [ "7^(2+1+2):GL2(7)" ],
    > [ "M11xA6.2^2" ],
    > [ "(S5xS5xS5):S3" ],
    > [ "(L2(11)xL2(11)):4" ],
    > [ "13^2:2.L2(13).4" ],
    > [ "(7^2:(3x2A4)xL2(7)).2" ],
    > [ "(13:6xL3(3)).2" ],
    > [ "13^(1+2):(3x4S4)" ],
    > [ "L2(71)" ],
    > [ "L2(59)" ],
    > [ "11^2:(5x2.A5)" ],
    > [ "L2(41)" ],
    > [ "L2(29).2" ],
    > [ "7^2:2psl(2,7)" ],
    > [ "L2(19).2" ],
    > [ "41:40" ],
    > ];;
\end{verbatim}

We compute upper bounds for the permutation character values
in the cases where the characters are not given explicitly.
(We could improve this by using additional information about the class
fusions, but this will not be necessary.)

\begin{verbatim}
    gap> for entry in primdata do
    >      s:= CharacterTable( entry[1] );
    >      if not IsBound( entry[2] ) then
    >        fus:= PossibleClassFusions( s, m );
    >        poss:= Set( List( fus, x -> InducedClassFunctionsByFusionMap( s, m,
    >                                      [ TrivialCharacter( s ) ], x )[1] ) );
    >        entry[2]:= List( [ 1 .. NrConjugacyClasses( m ) ],
    >                         i -> Maximum( List( poss, x -> x[i] ) ) );
    >      fi;
    >    od;
\end{verbatim}


According to~\cite{NW12},
any maximal subgroup of the Monster besides the above
{\num_known_monster_maxes} known classes
is an almost simple group whose socle is one of
L$_2(13)$, Sz$(8)$, U$_3(4)$, and U$_3(8)$.

We show that the elements of such subgroups are contained in
the union of $55$ conjugacy classes of the Monster
that cover less than one percent of the elements in the Monster.
For that, we compute the possible class fusions from the abovementioned
simple groups $S$ into the Monster, and then the possible class fusions
from the automorphic extensions of $S$ into the Monster,
using the possible class fusions of $S$.
(This approach is faster than computing each class fusion from scratch.)

After the following computations,
the list \verb|badclasses| will contain the positions of all those classes of $M$
that may contain elements in some of the hypothetical maximal subgroups.

For each simple group in question, we enter the identifiers of the
character tables of the automorphic extensions that can occur.
Note that the automorphism groups of the four groups have the structures
L$_2(13).2$, Sz$(8).3$, U$_3(4).4$, and U$_3(8).(3 \times S_3)$,
respectively.
We need not consider the groups U$_3(8).3^2$ and U$_3(8).(3 \times S_3)$
because already U$_3(8).3_2$ does not admit an embedding into $M$,
and we need not consider the group U$_3(8).S_3$ because its set of elements
is covered by its subgroups of the types U$_3(8).2$ and U$_3(8).3_2$.

\begin{verbatim}
    gap> PossibleClassFusions( CharacterTable( "U3(8).3_2" ), m );
    [  ]
    gap> badclasses:= [];;
    gap> names:= [
    >    [ "L2(13)", "L2(13).2" ],
    >    [ "Sz(8)", "Sz(8).3" ],
    >    [ "U3(4)", "U3(4).2", "U3(4).4" ],
    >    [ "U3(8)", "U3(8).2", "U3(8).3_1", "U3(8).3_2", "U3(8).3_3", "U3(8).6" ],
    >    ];;
    gap> for list in names do
    >      t:= CharacterTable( list[1] );
    >      tfusm:= PossibleClassFusions( t, m );
    >      UniteSet( badclasses, Flat( tfusm ) );
    >      for nam in list{ [ 2 .. Length( list ) ] } do
    >        ext:= CharacterTable( nam );
    >        for map1 in PossibleClassFusions( t, ext ) do
    >          inv:= InverseMap( map1 );
    >          for map2 in tfusm do
    >            init:= CompositionMaps( map2, inv );
    >            UniteSet( badclasses, Flat( PossibleClassFusions( ext, m,
    >                rec( fusionmap:= init ) ) ) );
    >          od;
    >        od;
    >      od;
    >    od;
    gap> badclasses;
    [ 1, 3, 4, 5, 6, 7, 9, 10, 11, 12, 14, 15, 17, 18, 19, 20, 21, 22, 24, 25, 
      27, 28, 30, 32, 33, 35, 36, 38, 39, 40, 42, 43, 44, 45, 46, 48, 49, 50, 51, 
      52, 53, 54, 55, 56, 60, 61, 62, 63, 70, 72, 73, 78, 82, 85, 86 ]
    gap> Length( badclasses );
    55
    gap> classlengths:= SizesConjugacyClasses( m );;
    gap> bad:= Sum( classlengths{ badclasses } ) / Size( m );;
    gap> Int( 10000 * bad ); 
    97
\end{verbatim}


\emph{In the original version of this file, also hypothetical maximal
subgroups with socle} L$_2(27)$ \emph{had been considered.
As a consequence, the list} \verb|badclasses| \emph{computed above had length $59$
in the original version;
the list contained also the classes at the positions $90, 94, 95$, and $96$,
that is, the classes} \verb|26B|, \verb|28B|, \verb|28C|, \verb|28D|.
\emph{The proportion} \verb|bad| \emph{of elements in the classes of $M$
described by}
\verb|badclasses| \emph{was about $2.05$ percent of $|M|$,
compared to the about $0.98$ percent in the current version.}


Now we estimate the lower bounds $\delta(s, g^G)$ introduced in
Section~\ref{characters}.
Let ${\cal{B}}$ denote the union of the classes described by \verb|badclasses|,
and let $\M$ denote a set of representatives of the
{\num_known_monster_maxes} known classes of maximal subgroups of $M$.

If $s \not\in {\cal{B}}$ then
\[
   \delta(s, g^G) =
     |s^G| - |s^G| \cdot \sum_{S \in \M} 1_S^M(s) \cdot 1_S^M(g) / 1_S^M(1) ,
\]
hence $\delta(s)$ can be computed from the corresponding primitive
permutation characters,
and a lower bound for $\delta(s)$ can be computed from the upper bounds for
the characters $1_S^G$ which are given by the list \verb|primdata|.

If $s \in {\cal{B}}$ then
the above equation for $\delta(s, g^G)$ holds at least for
$g \not\in {\cal{B}}$,
so $\sum_{g \in R \setminus {\cal{B}}} \delta(s, g^G)$ is a lower bound
for $\delta(s)$.
So \verb|primdata| yields a lower bound for $\delta(s)$ also for $s \in {\cal{B}}$,
by ignoring the pairs $(s, g)$ where both $s$ and $g$ lie in ${\cal{B}}$.

This means that modifying the output of \verb|LowerBoundsVertexDegrees|
as follows really yields lower bounds for the vertex degrees.
(Note that the row and column positions in the matrix returned by
\verb|LowerBoundsVertexDegrees| are shifted by one, compared to \verb|badclasses|.)

\begin{verbatim}
    gap> prim:= List( primdata, x -> x[2] );;
    gap> badpos:= Difference( badclasses, [ 1 ] ) - 1;;
    gap> bounds:= LowerBoundsVertexDegrees( classlengths, prim );;
    gap> for i in badpos do
    >      for j in badpos do
    >        bounds[i][j]:= 0;
    >      od;
    >    od;
\end{verbatim}

Now we sum up the bounds for the individual classes.
It turns out that the minimal vertex degree is more than $99$ percent
of $|M|$.
This proves that the generating graph of the Monster
satisfies P\'osa's criterion.

(And the minimal vertex degree of elements outside ${\cal{B}}$
is more than $99.99998$ percent of $|M|$.)

\emph{In the original version of this file,
we got only $97.95$ percent of $|M|$ as the lower bound
for the minimal vertex degree.
The bound for elements outside ${\cal{B}}$ was the same
in the original version.
The fact that the maximal subgroups of type} L$_2(41)$
\emph{had been ignored in the original version
did not affect the lower bound for the minimal vertex degree.}

\begin{verbatim}
    gap> degs:= List( bounds, Sum );;
    gap> Int( 10000 * Minimum( degs ) / Size( m ) );
    9902
    gap> goodpos:= Difference( [ 1 .. NrConjugacyClasses( m ) - 1 ], badpos );;
    gap> Int( 100000000 * Minimum( degs{ goodpos } ) / Size( m ) );
    99999987
\end{verbatim}

\subsection{Nonsimple Automorphism Groups of Sporadic Simple Groups}%
\label{sporaut}

Next we consider the nonsimple automorphism groups of the sporadic simple
groups.
Nontrivial outer automorphisms exist exactly in $12$ cases,
and then the simple group has index $2$ in its automorphism group.
The character tables of the groups and their maximal subgroups are
available in {\GAP}.

\begin{verbatim}
    gap> spornames:= AllCharacterTableNames( IsSporadicSimple, true,
    >                    IsDuplicateTable, false );;
    gap> sporautnames:= AllCharacterTableNames( IsSporadicSimple, true,
    >                       IsDuplicateTable, false,
    >                       OfThose, AutomorphismGroup );;
    gap> sporautnames:= Difference( sporautnames, spornames );
    [ "F3+.2", "Fi22.2", "HN.2", "HS.2", "He.2", "J2.2", "J3.2", "M12.2", 
      "M22.2", "McL.2", "ON.2", "Suz.2" ]
    gap> for tbl in List( sporautnames, CharacterTable ) do
    >      info:= HamiltonianCycleInfoFromCharacterTable( tbl );
    >      Print( Identifier( tbl ), ": ", info, "\n" );
    >    od;
    F3+.2: Chvatal for 0th closure, Posa for 1st closure
    Fi22.2: Chvatal for 0th closure, Posa for 1st closure
    HN.2: Chvatal for 0th closure, Posa for 1st closure
    HS.2: Chvatal for 1st closure, Posa for 2nd closure
    He.2: Chvatal for 0th closure, Posa for 1st closure
    J2.2: Chvatal for 0th closure, Posa for 1st closure
    J3.2: Chvatal for 0th closure, Posa for 1st closure
    M12.2: Chvatal for 0th closure, Posa for 1st closure
    M22.2: Posa for 1st closure
    McL.2: Chvatal for 0th closure, Posa for 1st closure
    ON.2: Chvatal for 0th closure, Posa for 1st closure
    Suz.2: Chvatal for 0th closure, Posa for 1st closure
\end{verbatim}

\subsection{Alternating and Symmetric Groups $A_n$, $S_n$,
for $5 \leq n \leq 13$}%
\label{symmalt}

For alternating and symmetric groups $A_n$ and $S_n$, respectively,
with $5 \leq n \leq 13$,
the table of marks or the character tables of the group
and all its maximal subgroups are available in {\GAP}.
So we can compute the character-theoretic bounds for vertex degrees.

\begin{verbatim}
    gap> for tbl in List( [ 5 .. 13 ], i -> CharacterTable(
    >                 Concatenation( "A", String( i ) ) ) )  do
    >      info:= HamiltonianCycleInfoFromCharacterTable( tbl );
    >      if info <> "Posa for 0th closure" then
    >        Print( Identifier( tbl ), ": ", info, "\n" );
    >      fi;
    >    od;
\end{verbatim}

No messages are printed, so the generating graphs of the alternating
groups in question satisfy P\'osa's criterion.

\begin{verbatim}
    gap> for tbl in List( [ 5 .. 13 ], i -> CharacterTable(
    >                 Concatenation( "S", String( i ) ) ) )  do
    >      info:= HamiltonianCycleInfoFromCharacterTable( tbl );
    >      Print( Identifier( tbl ), ": ", info, "\n" );
    >    od;
    A5.2: no decision
    A6.2_1: Chvatal for 4th closure, Posa for 5th closure
    A7.2: Posa for 1st closure
    A8.2: Chvatal for 2nd closure, Posa for 3rd closure
    A9.2: Chvatal for 2nd closure, Posa for 3rd closure
    A10.2: Chvatal for 2nd closure, Posa for 3rd closure
    A11.2: Posa for 1st closure
    A12.2: Chvatal for 2nd closure, Posa for 3rd closure
    A13.2: Posa for 1st closure
\end{verbatim}

We see that sufficiently large closures of the generating graphs of the
symmetric groups in question satisfy P\'osa's criterion,
except that the bounds for the symmetric group $S_5$ are not sufficient
for the proof.
In Section~\ref{smallalmostsimple}, it is shown that the 2nd closure of
the generating graph of $S_5$ satisfies P\'osa's criterion.

(We could find slightly better bounds derived only from character tables
which suffice to prove that the generating graph for $S_5$ contains a
Hamiltonian cycle, but this seems to be not worth while.)

\section{Computations With Groups}\label{grouptheor}

We prove in Section~\ref{smallsimp} that the generating graphs of the
nonabelian simple groups of order up to $10^6$ satisfy P\'osa's criterion,
and that the same holds for those nonabelian simple
groups of order between $10^6$ and $10^7$ that are not isomorphic
with some $\PSL(2,q)$.
(In Section~\ref{psl2q}, it is shown that the generating graph of
$\PSL(2,q)$ satifies P\'osa's criterion for any prime power $q$.)
Nonsimple nonsolvable groups of order up to $10^6$
are treated in Section~\ref{smallalmostsimple}.

(We could increase the bounds $10^6$ and $10^7$ with more computations
using the same methods.)

For our convenience, we provide a small function that takes as its
argument only the group in question,
and returns a string, the
return value of \verb|HamiltonianCycleInfo| for the vertex degrees computed from
the group.
(In order to speed up the computations,
the function computes the proper normal subgroups that contain the derived
subgroup of the given group, and enters the list of these groups as the
third argument of \verb|VertexDegreesGeneratingGraph|.)

\begin{verbatim}
    gap> HamiltonianCycleInfoFromGroup:= function( G )
    >     local ccl, nsg, der, degrees, classlengths;
    >     ccl:= ConjugacyClasses( G );
    >     if IsPerfect( G ) then
    >       nsg:= [];
    >     else
    >       der:= DerivedSubgroup( G );
    >       nsg:= Concatenation( [ der ],
    >                            IntermediateSubgroups( G, der ).subgroups );
    >     fi;
    >     degrees:= VertexDegreesGeneratingGraph( G, ccl, nsg );
    >     classlengths:= List( ccl, Size );
    >     return HamiltonianCycleInfo( classlengths, degrees );        
    > end;;
\end{verbatim}

\subsection{Nonabelian Simple Groups of Order up to $10^7$}\label{smallsimp}

Representatives of the $56$ isomorphism types of nonabelian simple groups
of order up to $10^6$ can be accessed in {\GAP} with the function
\verb|AllSmallNonabelianSimpleGroups|.

\begin{verbatim}
    gap> grps:= AllSmallNonabelianSimpleGroups( [ 1 .. 10^6 ] );;         
    gap> Length( grps );
    56
    gap> List( grps, StructureDescription );
    [ "A5", "PSL(3,2)", "A6", "PSL(2,8)", "PSL(2,11)", "PSL(2,13)", "PSL(2,17)", 
      "A7", "PSL(2,19)", "PSL(2,16)", "PSL(3,3)", "PSU(3,3)", "PSL(2,23)", 
      "PSL(2,25)", "M11", "PSL(2,27)", "PSL(2,29)", "PSL(2,31)", "A8", 
      "PSL(3,4)", "PSL(2,37)", "O(5,3)", "Sz(8)", "PSL(2,32)", "PSL(2,41)", 
      "PSL(2,43)", "PSL(2,47)", "PSL(2,49)", "PSU(3,4)", "PSL(2,53)", "M12", 
      "PSL(2,59)", "PSL(2,61)", "PSU(3,5)", "PSL(2,67)", "J1", "PSL(2,71)", "A9", 
      "PSL(2,73)", "PSL(2,79)", "PSL(2,64)", "PSL(2,81)", "PSL(2,83)", 
      "PSL(2,89)", "PSL(3,5)", "M22", "PSL(2,97)", "PSL(2,101)", "PSL(2,103)", 
      "HJ", "PSL(2,107)", "PSL(2,109)", "PSL(2,113)", "PSL(2,121)", "PSL(2,125)", 
      "O(5,4)" ]
    gap> for g in grps do                                             
    >      info:= HamiltonianCycleInfoFromGroup( g );
    >      if info <> "Posa for 0th closure" then
    >        Print( StructureDescription( g ), ": ", info, "\n" );
    >      fi;
    >    od;
\end{verbatim}

Nothing is printed during these computations,
so the generating graphs of all processed groups satisfy
P\'osa's criterion.

(On my notebook,
the above computations needed about 6300 seconds of CPU time.)

For simple groups of order larger than $10^6$, there is not such an easy
way (yet) to access representatives for each isomorphism type.
Therefore, first we compute the orders of nonabelian simple groups
between $10^6$ and $10^7$.

\begin{verbatim}
    gap> orders:= Filtered( [ 10^6+4, 10^6+8 .. 10^7 ],
    >                 n -> IsomorphismTypeInfoFiniteSimpleGroup( n ) <> fail );
    [ 1024128, 1123980, 1285608, 1342740, 1451520, 1653900, 1721400, 1814400, 
      1876896, 1934868, 2097024, 2165292, 2328648, 2413320, 2588772, 2867580, 
      2964780, 3265920, 3483840, 3594432, 3822588, 3940200, 4245696, 4680000, 
      4696860, 5515776, 5544672, 5663616, 5848428, 6004380, 6065280, 6324552, 
      6825840, 6998640, 7174332, 7906500, 8487168, 9095592, 9732420, 9951120, 
      9999360 ]
    gap> Length( orders );
    41
    gap> info:= List( orders, IsomorphismTypeInfoFiniteSimpleGroup );;
    gap> Number( info, x -> IsBound( x.series ) and x.series = "L"
    >                       and x.parameter[1] = 2 );
    31
\end{verbatim}

We see that there are $31$ groups of the type $\PSL(2,q)$ and $10$ other
nonabelian simple groups with order in the range from $10^6$ to $10^7$.
The former groups can be ignored because the generating graphs of any
group $\PSL(2,q)$ satisfies P\'osa's criterion,
see Section~\ref{psl2q}.
For the latter groups, we can apply the character-theoretic method
to prove that the generating graph satisfies P\'osa's criterion.

\begin{verbatim}
    gap> info:= Filtered( info, x -> not IsBound( x.series ) or
    >             x.series <> "L" or x.parameter[1] <> 2 );
    [ rec( name := "B(3,2) = O(7,2) ~ C(3,2) = S(6,2)", parameter := [ 3, 2 ],
          series := "B" ), rec( name := "A(10)", parameter := 10, series := "A" ),
      rec( name := "A(2,7) = L(3,7) ", parameter := [ 3, 7 ], series := "L" ),
      rec( name := "2A(3,3) = U(4,3) ~ 2D(3,3) = O-(6,3)", parameter := [ 3, 3 ],
          series := "2A" ), rec( name := "G(2,3)", parameter := 3, series := "G" )
        ,
      rec( name := "B(2,5) = O(5,5) ~ C(2,5) = S(4,5)", parameter := [ 2, 5 ],
          series := "B" ),
      rec( name := "2A(2,8) = U(3,8)", parameter := [ 2, 8 ], series := "2A" ),
      rec( name := "2A(2,7) = U(3,7)", parameter := [ 2, 7 ], series := "2A" ),
      rec( name := "A(3,3) = L(4,3) ~ D(3,3) = O+(6,3) ", parameter := [ 4, 3 ],
          series := "L" ),
      rec( name := "A(4,2) = L(5,2) ", parameter := [ 5, 2 ], series := "L" ) ]
    gap> names:= [ "S6(2)", "A10", "L3(7)", "U4(3)", "G2(3)", "S4(5)", "U3(8)",
    >              "U3(7)", "L4(3)", "L5(2)" ];;
    gap> for tbl in List( names, CharacterTable ) do
    >      info:= HamiltonianCycleInfoFromCharacterTable( tbl );
    >      if info <> "Posa for 0th closure" then
    >        Print( Identifier( tbl ), ": ", info, "\n" );
    >      fi;
    >    od;
\end{verbatim}

\subsection{Nonsimple Groups with Nonsolvable Socle of Order at most $10^6$}%
\label{smallalmostsimple}

Let $G$ be a nonsolvable group such that $G/N$ is cyclic for all nontrivial
normal subgroups $N$ of $G$.
Then the socle Soc$(G)$ of $G$ is the unique minimal normal subgroup.
Moreover, Soc$(G)$ is nonsolvable and thus a direct product of isomorphic
nonabelian simple groups,
and $G$ is isomorphic to a subgroup of Aut$($Soc$(G))$.

In order to deal with all such groups $G$ for which additionally
$|$Soc$(G)| \leq 10^6$ holds,
it is sufficient to run over the simple groups $S$ of order up to $10^6$
and to consider those subgroups $G$ of Aut$(S^n)$,
with $|S|^n \leq 10^6$, for which Inn$(G)$ is the unique minimal normal
subgroups and $G / $Inn$(G)$ is cyclic.

We show that for each such group,
a sufficient closure of the generating graph satisfies P\'osa's criterion.

\begin{verbatim}
    gap> grps:= AllSmallNonabelianSimpleGroups( [ 1 .. 10^6 ] );;         
    gap> for simple in grps do
    >      for n in [ 1 .. LogInt( 10^6, Size( simple ) ) ] do
    >        # Compute the n-fold direct product S^n.
    >        soc:= CallFuncList( DirectProduct,
    >                            ListWithIdenticalEntries( n, simple ) );
    >        # Compute Aut(S^n) as a permutation group.
    >        aut:= Image( IsomorphismPermGroup( AutomorphismGroup( soc ) ) );
    >        aut:= Image( SmallerDegreePermutationRepresentation( aut ) );
    >        # Compute class representatives of subgroups of Aut(S^n)/Inn(S^n).
    >        socle:= Socle( aut );
    >        epi:= NaturalHomomorphismByNormalSubgroup( aut, socle );
    >        # Compute the candidates for G.
    >        # (By the above computations, we need not consider simple groups.)
    >        reps:= List( ConjugacyClassesSubgroups( Image( epi ) ),
    >                     Representative );
    >        reps:= Filtered( reps, x -> IsCyclic( x ) and Size( x ) <> 1 );
    >        greps:= Filtered( List( reps, x -> PreImages( epi, x ) ),
    >                          x -> Length( MinimalNormalSubgroups( x ) ) = 1 );
    >        for g in greps do
    >          # We have to deal with a *transitive* permutation group.
    >          # (Each group in question acts faithfully on an orbit.)
    >          if not IsTransitive( g ) then
    >            g:= First( List( Orbits( g, MovedPoints( g ) ),
    >                             x -> Action( g, x ) ),
    >                       x -> Size( x ) = Size( g ) );
    >          fi;
    >          # Check this group G.
    >          info:= HamiltonianCycleInfoFromGroup( g );
    >          Print( Name( simple ), "^", n, ".", Size( g ) / Size( soc ), ": ",
    >                 info, "\n" );
    >        od;
    >      od;
    >    od;
    A5^1.2: Posa for 2nd closure
    A5^2.2: Posa for 0th closure
    A5^2.4: Posa for 0th closure
    A5^3.3: Posa for 0th closure
    A5^3.6: Chvatal for 1st closure, Posa for 2nd closure
    PSL(2,7)^1.2: Chvatal for 0th closure, Posa for 1st closure
    PSL(2,7)^2.2: Posa for 0th closure
    PSL(2,7)^2.4: Posa for 0th closure
    A6^1.2: Chvatal for 0th closure, Posa for 1st closure
    A6^1.2: Chvatal for 4th closure, Posa for 5th closure
    A6^1.2: Chvatal for 0th closure, Posa for 1st closure
    A6^2.2: Posa for 0th closure
    A6^2.4: Posa for 0th closure
    A6^2.4: Posa for 0th closure
    A6^2.4: Posa for 0th closure
    PSL(2,8)^1.3: Posa for 0th closure
    PSL(2,8)^2.2: Posa for 0th closure
    PSL(2,8)^2.6: Chvatal for 0th closure, Posa for 1st closure
    PSL(2,11)^1.2: Chvatal for 0th closure, Posa for 1st closure
    PSL(2,11)^2.2: Posa for 0th closure
    PSL(2,11)^2.4: Posa for 0th closure
    PSL(2,13)^1.2: Chvatal for 0th closure, Posa for 1st closure
    PSL(2,17)^1.2: Chvatal for 0th closure, Posa for 1st closure
    A7^1.2: Posa for 1st closure
    PSL(2,19)^1.2: Chvatal for 0th closure, Posa for 1st closure
    PSL(2,16)^1.2: Chvatal for 0th closure, Posa for 1st closure
    PSL(2,16)^1.4: Chvatal for 0th closure, Posa for 1st closure
    PSL(3,3)^1.2: Chvatal for 0th closure, Posa for 1st closure
    PSU(3,3)^1.2: Chvatal for 0th closure, Posa for 1st closure
    PSL(2,23)^1.2: Chvatal for 0th closure, Posa for 1st closure
    PSL(2,25)^1.2: Chvatal for 0th closure, Posa for 1st closure
    PSL(2,25)^1.2: Chvatal for 0th closure, Posa for 1st closure
    PSL(2,25)^1.2: Chvatal for 0th closure, Posa for 1st closure
    PSL(2,27)^1.2: Chvatal for 0th closure, Posa for 1st closure
    PSL(2,27)^1.3: Posa for 0th closure
    PSL(2,27)^1.6: Chvatal for 0th closure, Posa for 1st closure
    PSL(2,29)^1.2: Chvatal for 0th closure, Posa for 1st closure
    PSL(2,31)^1.2: Chvatal for 0th closure, Posa for 1st closure
    A8^1.2: Chvatal for 2nd closure, Posa for 3rd closure
    PSL(3,4)^1.2: Chvatal for 0th closure, Posa for 1st closure
    PSL(3,4)^1.2: Chvatal for 1st closure, Posa for 2nd closure
    PSL(3,4)^1.2: Chvatal for 0th closure, Posa for 1st closure
    PSL(3,4)^1.3: Posa for 0th closure
    PSL(3,4)^1.6: Chvatal for 0th closure, Posa for 1st closure
    PSL(2,37)^1.2: Chvatal for 0th closure, Posa for 1st closure
    PSp(4,3)^1.2: Chvatal for 1st closure, Posa for 2nd closure
    Sz(8)^1.3: Posa for 0th closure
    PSL(2,32)^1.5: Posa for 0th closure
    PSL(2,41)^1.2: Chvatal for 0th closure, Posa for 1st closure
    PSL(2,43)^1.2: Chvatal for 0th closure, Posa for 1st closure
    PSL(2,47)^1.2: Chvatal for 0th closure, Posa for 1st closure
    PSL(2,49)^1.2: Chvatal for 0th closure, Posa for 1st closure
    PSL(2,49)^1.2: Chvatal for 0th closure, Posa for 1st closure
    PSL(2,49)^1.2: Chvatal for 0th closure, Posa for 1st closure
    PSU(3,4)^1.2: Chvatal for 0th closure, Posa for 1st closure
    PSU(3,4)^1.4: Chvatal for 0th closure, Posa for 1st closure
    PSL(2,53)^1.2: Chvatal for 0th closure, Posa for 1st closure
    M12^1.2: Chvatal for 0th closure, Posa for 1st closure
    PSL(2,59)^1.2: Chvatal for 0th closure, Posa for 1st closure
    PSL(2,61)^1.2: Chvatal for 0th closure, Posa for 1st closure
    PSU(3,5)^1.2: Chvatal for 0th closure, Posa for 1st closure
    PSU(3,5)^1.3: Posa for 0th closure
    PSL(2,67)^1.2: Chvatal for 0th closure, Posa for 1st closure
    PSL(2,71)^1.2: Chvatal for 0th closure, Posa for 1st closure
    A9^1.2: Chvatal for 2nd closure, Posa for 3rd closure
    PSL(2,73)^1.2: Chvatal for 0th closure, Posa for 1st closure
    PSL(2,79)^1.2: Chvatal for 0th closure, Posa for 1st closure
    PSL(2,64)^1.2: Chvatal for 0th closure, Posa for 1st closure
    PSL(2,64)^1.3: Posa for 0th closure
    PSL(2,64)^1.6: Chvatal for 0th closure, Posa for 1st closure
    PSL(2,81)^1.2: Chvatal for 0th closure, Posa for 1st closure
    PSL(2,81)^1.2: Chvatal for 0th closure, Posa for 1st closure
    PSL(2,81)^1.2: Chvatal for 0th closure, Posa for 1st closure
    PSL(2,81)^1.4: Chvatal for 0th closure, Posa for 1st closure
    PSL(2,81)^1.4: Chvatal for 0th closure, Posa for 1st closure
    PSL(2,83)^1.2: Chvatal for 0th closure, Posa for 1st closure
    PSL(2,89)^1.2: Chvatal for 0th closure, Posa for 1st closure
    PSL(3,5)^1.2: Chvatal for 0th closure, Posa for 1st closure
    M22^1.2: Posa for 1st closure
    PSL(2,97)^1.2: Chvatal for 0th closure, Posa for 1st closure
    PSL(2,101)^1.2: Chvatal for 0th closure, Posa for 1st closure
    PSL(2,103)^1.2: Chvatal for 0th closure, Posa for 1st closure
    J_2^1.2: Chvatal for 0th closure, Posa for 1st closure
    PSL(2,107)^1.2: Chvatal for 0th closure, Posa for 1st closure
    PSL(2,109)^1.2: Chvatal for 0th closure, Posa for 1st closure
    PSL(2,113)^1.2: Chvatal for 0th closure, Posa for 1st closure
    PSL(2,121)^1.2: Chvatal for 0th closure, Posa for 1st closure
    PSL(2,121)^1.2: Chvatal for 0th closure, Posa for 1st closure
    PSL(2,121)^1.2: Chvatal for 0th closure, Posa for 1st closure
    PSL(2,125)^1.2: Chvatal for 0th closure, Posa for 1st closure
    PSL(2,125)^1.3: Posa for 0th closure
    PSL(2,125)^1.6: Chvatal for 0th closure, Posa for 1st closure
    PSp(4,4)^1.2: Chvatal for 0th closure, Posa for 1st closure
    PSp(4,4)^1.4: Posa for 0th closure
\end{verbatim}

\section{The Groups $\PSL(2,q)$}\label{psl2q}

We show that the generating graph of any group $\PSL(2,q)$, for $q \geq 2$,
satisfies P\'osa's criterion.
Throughout this section,
let $q = p^f$ for a prime integer $p$, and $G = \PSL(2,q)$.
Set $k = \gcd(q-1, 2)$.

\begin{lem}[{{see~\cite[II., \S~8]{Hup67}}}]
The subgroups of $G$ are
\begin{itemize}
\item[(1)]
    cyclic groups of order dividing $(q \pm 1)/k$,
    and their normalizers, which are dihedral groups of order $2 (q \pm 1)/k$,
\item[(2)]
    subgroups of Sylow $p$ normalizers, which are semidirect products of
    elementary abelian groups of order $q$ with cyclic groups of order
    $(q-1)/k$,
\item[(3)]
    subgroups isomorphic with $\PSL(2, p^m)$ if $m$ divides $f$,
    and isomorphic with $\PGL(2, p^m)$ if $2 m$ divides $f$,
\item[(4)]
    subgroups isomorphic with $A_4$, $S_4$, or $A_5$,
    for appropriate values of $q$.
\end{itemize}
$G$ contains exactly one conjugacy class of cyclic subgroups of each of
the orders $(q-1)/k$ and $(q+1)/k$,
and each nonidentity element of $G$ is contained in exactly one of these
subgroups or in exactly one Sylow $p$ subgroup of $G$.
\end{lem}

We estimate the number of elements that are contained in subgroups
of type~(3).

\begin{lem}
Let $n_{sf}(q)$ denote the number of those nonidentity elements in
$G$ that are contained in proper subgroups of type~(3).
Then $n_{sf}(q) \leq q^2 (\frac{2 p}{p-1} (\sqrt{q}-1) - 1)$.
If $f$ is a prime then $n_{sf}(q) \leq (2p-1) q^2$ holds,
and if $p = q$ then we have of course $n_{sf}(q) = 0$.
\end{lem}

\tthdump{\begin{proof}}
The group $\PGL(2, p^m)$ is equal to $\PSL(2, p^m)$ for $p = 2$,
and contains $\PSL(2, p^m)$ as a subgroup of index two if $p \not= 2$.
So the largest element order in $\PGL(2, p^m)$ is at most $p^m+1$.
Let $C$ be a cyclic subgroup of order $(q + \epsilon)/k$ in $G$,
for $\epsilon \in \{ \pm 1 \}$.
The intersection of $C$ with any subgroup of $G$ isomorphic with
$\PGL(2, p^m)$ or $\PSL(2, p^m)$
is contained in the union of the unique subgroups of the orders
$\gcd(|C|, p^m + 1)$ and $\gcd(|C|, p^m - 1)$ in $C$.
So $C$ contains at most $2 p^m - 2$ nonidentity elements that can lie inside
subgroups isomorphic with $\PGL(2, p^m)$ or $\PSL(2, p^m)$.
Hence $C$ contains at most $\sum_m (2 p^m - 2)$ nonidentity elements
in proper subgroups of type~(3),
where $m$ runs over the proper divisors of $f$.
This sum is bounded from above by
$\sum_{m=1}^{f/2} (2 p^m - 2) \leq \frac{2 p}{p-1} (\sqrt{q}-1) - 2$.

The numbers of cyclic subgroups of the orders $(q + \epsilon)/k$ in $G$
are $q (q - \epsilon) / 2$, so $G$ contains altogether $q^2$ such cyclic
subgroups.
They contain at most $q^2 (\frac{2 p}{p-1} (\sqrt{q}-1) - 2)$
elements inside proper subgroups of the type (3).

All elements of order $p$ in $G$ are contained in subgroups of type~(3),
and there are exactly $q^2 - 1$ such elements.
This yields the claimed bound for $n_{sf}(q)$.
The better bound for the case that $f$ is a prime follows from
$\sum_m (2 p^m - 2) = 2 p - 2$ if $m$ ranges over the proper divisors of $f$.
\tthdump{\end{proof}}

Using these bounds, we see that the vertex degree of any element in $G$
that does not lie in subgroups of type~(4) is larger than $|G|/2$.
(In fact we could use the calculations below to derive a better
asymptotic bound, but this is not an issue here.)

\begin{lem}
Let $s \in G$ be an element of order larger than $5$.
Then $|\{ g \in G; \langle g, s \rangle = G \}| > |G|/2$.
\end{lem}

\tthdump{\begin{proof}}
First suppose that the order of $s$ divides $(q+1)/k$ or $(q-1)/k$.
If $g \in G$ such that $U = \langle s, g \rangle$
is a proper subgroup of $G$ then $U \leq N_G(\langle s \rangle)$ or
$U$ lies in a Sylow $p$ normalizer of $G$ or $U$ lies in a subgroup
of type~(3).
Since $s$ is contained in at most two Sylow $p$ normalizers
(each Sylow $p$ normalizer contains $q$ cyclic subgroups of order $(q-1)/k$,
and $G$ contains $q+1$ Sylow normalizers and $q (q+1)/2$ cyclic subgroups
of order $(q-1)/k$),
the number of $g \in G$ with the property that
$\langle s, g \rangle \not= G$ is at most
$N = 2(q+1)/k + 2 q(q-1)/k + n_{sf}(q) = 2(q^2+1)/k + n_{sf}(q)$;
for $s$ of order equal to $(q+1)/k$ or $(q-1)/k$, we can set
$N = 2(q^2+1)/k$.

Any element $s$ of order $p$ (larger than $5$), lies only in a unique
Sylow $p$ normalizer and in subgroups of type~(3), so the bound $N$
holds also in this case.

For $f = 1$, $N$ is smaller than $|G|/2 = q (q^2-1) / (2 k)$
if $q \geq 5$.
(The statement of the lemma is trivially true for $q \leq 5$.)

For primes $f$, $N$ is smaller than $|G|/2$
if $q^2 (q-8p) > q+4$ holds, which is true for $p^f > 8p$.
Only the following values of $p^f$ with prime $f$ do not satisfy this
condition:
$2^2$ and $3^2$ (where no element of order larger than $5$ exists),
$2^3$ (where only elements of order equal to $q \pm 1$ must be considered),
$5^2$ and $7^2$ (where $n_{sf}(q) < (p-1) q (q+1)$ because
in these cases the cyclic subgroups of order $(q+1)/k$ cannot contain
nonidentity elements in subgroups of type~(3)).

Finally, if $f$ is not a prime then $N$ is smaller than $|G|/2$
if $q^2 (q - \frac{8p}{p-1} (\sqrt{q}-1)) > q+4$ holds, which is true
for $q \geq 256$.
The only values of $p^f$ with non-prime $f$ that do not satisfy this condition
are $2^4$, $2^6$, and $3^4$.
In all three cases, we have in fact $N < |G|/2$,
where we have to use the better bound $n_{sf}(q) < 16 q^2$ in the third case.
\tthdump{\end{proof}}

In order to show that the generating graph of $G$ satisfies
P\'osa's criterion,
it suffices to show that the vertex degrees of involutions is larger than
the number of involutions, and that the vertex degrees of elements of orders
$2$, $3$, $4$, and $5$ are larger than the number of elements whose order
is at most $5$.

\begin{lem}
Let $n(q, m)$ denote the number of elements of order $m$ in $G$,
and let $\varphi(m)$ denote the number of prime residues modulo $m$.
\begin{itemize}
\item
    We have $n(q, 2) = q^2 - 1$ if $q$ is even
    and $n(q, 2) \leq q (q+1)/2$ if $q$ is odd.
\item
    For $m \in \{ 3, 4, 5 \}$, we have $n(q, m) \leq \varphi(m) q (q+1)/2$.
\item
    We have $n(q, (q+1)/k) = \varphi((q+1)/k) q (q-1)/2$.
\end{itemize}
\end{lem}


\begin{lem}\label{order2345}
If $q > 11$
then each involution in $G$ has vertex degree larger than $n(q, 2)$.

If $\varphi((q+1)/k) \geq 12$ then each element of order $3$, $4$, or $5$
has vertex degree larger than $\sum_{m=2}^5 n(q, m)$.
\end{lem}

\tthdump{\begin{proof}}
Let $s \in G$ of order at most $5$.
For each element $g \in G$ of order $(q+1)/k$,
$U = \langle g, s \rangle$ is either $G$ or contained in the dihedral group
of order $2(q+1)/k$ that normalizes $\langle g \rangle$.

If $s$ is an involution then
the number of such dihedral groups that contain $s$ is at most $(q+3)/2$,
and at least
$n(q, (q+1)/k) - \varphi((q+1)/k) (q+3)/2 = \varphi((q+1)/k) (q^2-2q-3)/2$
elements of order $(q+1)/k$ contribute to the vertex degree of $s$.
This number is larger than $q^2 - 1 \geq n(q, 2)$
if $q > 11$ (and hence $\varphi((q+1)/k) \geq 3$) holds.

If $s$ is an element of order $3$, $4$, or $5$
then $U \not= G$ means that $s \in \langle g \rangle$,
so at least
$n(q, (q+1)/k) - 4$ elements of order $(q+1)/k$
contribute to the vertex degree of $s$.
This number is larger than $5 q (q+1) > \sum_{m=2}^5 n(q, m)$
if $\varphi((q+1)/k) \geq 12$.
\tthdump{\end{proof}}

It remains to deal with the values $q$ where $\varphi((q+1)/k) < 12$,
that is, $(q+1)/k \leq 30$.
We compute that the statement of Lemma~\ref{order2345} is true also for
prime powers $q$ with $11 < q \leq 59$.

\begin{verbatim}
    gap> TestL2q:= function( t )
    >    local name, orders, nccl, cl, prim, bds, n, ord;
    > 
    >    name:= Identifier( t );
    >    orders:= OrdersClassRepresentatives( t );
    >    nccl:= Length( orders );
    >    cl:= SizesConjugacyClasses( t );
    >    prim:= PrimitivePermutationCharacters( t );
    >    bds:= List( LowerBoundsVertexDegrees( cl, prim ), Sum );
    >    n:= List( [ 1 .. 5 ], i -> Sum( cl{ Filtered( [ 1 .. nccl ],
    >                                        x -> orders[x] = i ) } ) );
    >    if ForAny( Filtered( [ 1 .. nccl ], i -> orders[i] > 5 ),
    >               i -> bds[i-1] <= Size( t ) / 2 ) then
    >      Error( "problem with large orders for ", name );
    >    elif ForAny( Filtered( [ 1 .. nccl ], i -> orders[i] = 2 ),
    >                 i -> bds[i-1] <= n[2] ) then
    >      Error( "problem with order 2 for ", name, "\n" );
    >    elif ForAny( Filtered( [ 1 .. nccl ], i -> orders[i] in [ 3 .. 5 ] ),
    >                 i -> bds[i-1] <= Sum( n{ [ 2 .. 5 ] } ) ) then
    >      Error( "problem with order in [ 3 .. 5 ] for ", name );
    >    fi;
    > end;;
    gap> for q in Filtered( [ 13 .. 59 ], IsPrimePowerInt ) do
    >      TestL2q( CharacterTable( Concatenation( "L2(", String( q ), ")" ) ) );
    >    od;
\end{verbatim}

For $2 \leq q \leq 11$, the statement of Lemma~\ref{order2345} is not true
but P\'osa's criterion is satisfied for the generating graphs of the
groups $\PSL(2,q)$ with $2 \leq q \leq 11$.

\begin{verbatim}
    gap> for q in Filtered( [ 2 .. 11 ], IsPrimePowerInt ) do
    >      info:= HamiltonianCycleInfoFromGroup( PSL( 2, q ) );
    >      if info <> "Posa for 0th closure" then
    >        Print( q, ": ", info, "\n" );
    >      fi;
    >    od;
\end{verbatim}

\tthdump{\addcontentsline{toc}{section}{References}}

\bibliographystyle{amsalpha}
\newcommand{\etalchar}[1]{$^{#1}$}
\providecommand{\bysame}{\leavevmode\hbox to3em{\hrulefill}\thinspace}
\providecommand{\MR}{\relax\ifhmode\unskip\space\fi MR }
\providecommand{\MRhref}[2]{%
  \href{http://www.ams.org/mathscinet-getitem?mr=#1}{#2}
}
\providecommand{\href}[2]{#2}


\end{document}